\documentclass[12pt]{amsart}

\usepackage{amsmath,amssymb,amscd,amsthm}
\usepackage{latexsym}
\usepackage[dvips]{graphics,epsfig}
\newtheorem{theorem}{Theorem}
\newtheorem{lemma}{Lemma}
\newtheorem{proposition}{Proposition}
\newtheorem{remark}{Remark}

\newtheorem{corollary}{Corollary}

\renewcommand{\le}{\leqslant}
\renewcommand{\ge}{\geqslant}

\newcommand{\vol}{{\rm vol}}
\newcommand{\af}{{\rm af}}

\newcommand{\rtwo}{\mathbb R^2}

\newcommand{\rn}{\mathbb R^n}
\newcommand{\Prod}{\mathcal P}

\newcommand{\beq}{\begin{equation}}
\newcommand{\eeq}{\end{equation}}
\newcommand{\beqna}{\begin{eqnarray*}}
\newcommand{\eeqna}{\end{eqnarray*}}
\newcommand{\beqn}{\begin{equation*}}
\newcommand{\eeqn}{\end{equation*}}
\newcommand{\bp}{\begin{proof}}
\newcommand{\ep}{\end{proof}}
\newcommand{\bprop}{\begin{proposition}}
\newcommand{\eprop}{\end{proposition}}
\newcommand{\bt}{\begin{theorem}}
\newcommand{\et}{\end{theorem}}
\newcommand{\bex}{\begin{Example}}
\newcommand{\eex}{\end{Example}}
\newcommand{\bc}{\begin{corollary}}
\newcommand{\ec}{\end{corollary}}
\newcommand{\bl}{\begin{lemma}}
\newcommand{\el}{\end{lemma}}
\newcommand{\br}{\begin{remark}}
\newcommand{\er}{\end{remark}}

\begin{document}

\title
[Mahler conjecture: local minimality of the unit cube] {A remark on
the Mahler conjecture: local minimality of the unit cube}

\author{}

\author[F. Nazarov, F. Petrov, D. Ryabogin, and A. Zvavitch]{Fedor Nazarov, Fedor Petrov, Dmitry Ryabogin, and Artem Zvavitch}

\address{Department of Mathematics,
University of Wisconsin, Madison 480 Lincoln Drive Madison, WI 53706
}\email{nazarov@math.wisc.edu}

\address{St. Petersburg Department of Steklov Institute of
Mathematics} \email{fedyapetrov@gmail.com}
\address{Department of Mathematics, Kent State University,
Kent, OH 44242, USA} \email{ryabogin@math.kent.edu}
\address{Department of Mathematics, Kent State University,
Kent, OH 44242, USA} \email{zvavitch@math.kent.edu}
\thanks{Supported in part by U.S.~National Science Foundation grants
 DMS-0652684, DMS-0800243, DMS-0808908.} 

\date{\today}

\begin{abstract}
We prove  that  the unit cube $B^n_{\infty}$ is a strict local minimizer for
the Mahler volume product $vol_n(K)vol_n(K^*)$  in the class of origin symmetric convex bodies endowed with the Banach-Mazur distance. \end{abstract}

\maketitle

\section{Introduction}
In 1939 Mahler \cite{Ma}  asked  the following question.
 Let
 $K \subset \rn$, $n\ge 2$,  be a convex origin-symmetric body and let
$$
K^*:= \{ \xi \in \rn: x \cdot \xi \le  1 \,\,\forall x \in K \}
$$
be its polar body.  Define $\Prod(K)=\vol_n(K)\vol_n(K^*)$. Is it true that we always have
$$
 \Prod(K)\, \ge \Prod(B^n_{\infty}),
$$
where  $B_\infty^n=\{x\in \rn: |x_i| \le 1, 1 \le i \le n\}$?

Mahler himself proved in \cite{Ma} that the answer is affirmative when $n=2$.
 There are
several other proofs of the two-dimensional result, see for example
the proof of M. Meyer, \cite{Me2}, but the question  is still open
even in the three-dimensional case.

In the $n$-dimensional case, the conjecture  has been  verified for
some special classes of bodies,  namely, for bodies  that are unit
balls of Banach spaces with $1$-unconditional bases, \cite{SR},
\cite{R2}, \cite{Me1}, and for zonoids, \cite{R1}, \cite{GMR}.

Bourgain and Milman \cite{BM} proved the inequality
$$
 \Prod(K)^{1/n} \ge c \Prod(B^n_{\infty})^{1/n},
$$
 with some constant $c>0$ independent of $n$. The best known constant $c=\pi/4$ is due to Kuperberg \cite{Ku}.

Note that the exact upper bound for $\Prod(K)$ is known:
$$
\Prod(K)\le \Prod(B_2^n),
$$
where $B_2^n$ is the $n$-dimensional Euclidean unit  ball. This
bound was proved by Santalo \cite{Sa}. In  \cite{MeP}  it was shown
that the equality holds only if $K$ is an ellipsoid.

Let $d_{BM}(K,L)=\inf\{b/a:\,\exists T \in GL(n) \mbox{  such that
 } aK\subseteq TL\subseteq b K\}$ be the Banach-Mazur
multiplicative distance between bodies $K, L \subset \rn$.
 In this paper we prove the
following result.

\noindent{\bf Theorem.}   {\it
Let $K\subset \rn$ be an origin-symmetric convex body. Then
$$
\Prod(K) \ge \Prod(B^n_{\infty}) ,
$$
provided that $d_{BM}(K,B^n_{\infty})\le 1+\delta$,  and $\delta=\delta(n)>0$
is small enough. Moreover, the equality holds only if $d_{BM}(K, B_\infty^n)=1$, i.e.,  if $K$ is a parallelepiped.}\\

 {\bf Acknowledgment}. We are indebted to Matthew Meyer for
valuable discussions.\\

{\bf Notation}. Given a set $F\subset \rn$, we define $\af(F)$ to be
the affine subspace of the minimal dimension containing $F$, and
$l(F)$ to be the  linear subspace parallel to $\af(F)$ of the same
dimension. The boundary of a convex body $K$ is denoted by $\partial
K$. For a given set $P\subset \rn$, we write $P^{\perp}=\{x\in \rn
:\,x\cdot y=0,\forall y\in P\}$. Let ${\mathcal F}$ be the set of
all faces $F$ of all dimensions of the cube $B_\infty^n$. We  denote
by $c_F$ the center of a face $F \in {\mathcal F}$. We also denote
$B_p^n=\{x\in \rn: \sum_i |x_i|^p \le 1 \}$. By $C$ and $c$ we denote
large and small positive constants that may change from line to line and may
depend on the dimension $n$.

\section{Auxiliary results}

Note that $\Prod(TK)=\Prod(K)$ for all $T\in GL(n)$. We will use  this fact for choosing a canonical position for $K$.

 \bl\label{l:tangent}  Let $P$ be a parallelepiped of minimal volume
 containing a
  convex origin-symmetric body $K$.  Let $T:\rn\to\rn$
be a linear transformation such that $P=TB_\infty^n$. Then  $T^{-1}K
\subset B_\infty^n$ and $\pm e_j\in \partial T^{-1}K$, $j=1,...,n$.
\el

\bp Note that $B_\infty^n$ is a parallelepiped of minimal volume
containing $T^{-1}K$. If $e_j \not \in T^{-1}K$, then there exists an
affine hyperplane $H\ni e_j$ such that $H\cap T^{-1} K=\varnothing$.
Note that the volume of the parallelepiped bounded by $H, -H$, and
the affine hyperplanes $\{x: x\cdot e_i=\pm 1\}$, $i\not =j$, equals
$\vol_n(B_\infty^n)$, and that this parallelepiped still contains
$K$. But then we can shift $H$ and $-H$ towards $K$ a little bit and
a get a new parallelepiped of smaller volume containing K.
 \ep

 We shall need the following simple technical lemma.

 \bl\label{l:small}  Let $P\subset \rn$ be a star-shaped (with respect to the origin) polytope such that every its $(n-1)$-dimensional face $F$ has area at least $A$ and satisfies
  $d(\af(F),0) \ge r$. Let $x\not \in (1+\delta)P$ for some $\delta >0$. Then
 $$
 \vol_n({\rm conv}(P,x)) \ge \vol_n (P) +\frac{\delta r A}{n}.
 $$
 \el
 \bp Let $y=\partial P \cap [0,x]$. Let $F$ be a face of $P$ containing $y$. Then ${\rm conv}(P,x)\setminus P$ contains the pyramid with base $F$ and apex $x$.
 The assumptions of the lemma imply that the height of this pyramid is at least $\delta \, d(\af(F), 0) \ge \delta r$, so its volume is at least $\frac{\delta r A}{n}$.
\ep

If $K$ is sufficiently close to $B^n_{\infty}$, then  $K$ is also close to the parallelepiped of minimal volume
containing $K$.

 \bl\label{l:min}  Let  $K$ be a convex body
satisfying
$$
(1-\delta)B^n_{\infty}\subset K \subset  B^n_{\infty}.
$$
Then there exists a constant $C$ and a linear operator $T$ such
that
$$
(1-C \delta)B^n_{\infty}\subset T^{-1} K \subset  B^n_{\infty},
$$
 and $\pm e_i \in T^{-1}
K$. \el

\bp  Let as before $P=T B_\infty^n$ be a parallelepiped of minimal
volume containing $K$.  Note that $\vol_n(P) \le 2^n$. On the other
hand, if $x\in P\setminus (1+\kappa)(1-\delta)B_\infty^n$, then, by
Lemma \ref{l:small},
$$
\vol_n(P) \ge 2^n(1-\delta)^n + \kappa \frac{2^{n-1}}{n}
(1-\delta)^n.
$$
The right hand side is greater than $2^n$ if $\kappa> \kappa_0= 2n ((1-\delta)^{-n}-1)$. Thus $P\subset (1+\kappa_0)(1-\delta)B_\infty^n$ and thereby $(1-\kappa_0)P \subset (1-\delta)B_\infty^n \subset K$. It remains to note that $\kappa_0 \le 4n^2 \delta$ for sufficiently small $\delta>0$.
 \ep

Thus, replacing $K$ by  its suitable linear image we may assume
everywhere below that $K\subset B_\infty^n$, $\pm e_j\in \partial
K$, $j=1, \dots, n$. Let $\delta>0$ be the minimal number such that
$(1-\delta) B_\infty^n \subset K$.

\section{Computation of  the kernel of the  differential of the volume function}

 Choose some numbers $a_k >0$, $k=0, \dots, n-1$, and define the polytope $Q_0$ as the union of the  simplices
$$
S_{\mathbb F}={\rm conv} (0, a_0 c_{F_0}, a_1 c_{F_1}, \dots, a_{n-1} c_{F_{n-1}}),
$$
where  ${\mathbb F}= \{F_0,\dots, F_{n-1}\}$  runs over all flags ($F_0\subset F_1\subset F_2 \subset \dots\subset F_{n-1}, {\rm dim} F_j=j$) of faces of the unit cube.

Choose now some points $x_F$ close to $x_F^0=a_{{\rm dim} F} c_F$
and consider the polytope  $Q$ defined in the same way using the
points $x_F$. Consider the function $g(\{x_F\}_{F\in {\mathcal
F}})=\vol_n(Q)$. It is just a polynomial of degree $n$ of the
coordinates of $x_F$, so it is infinitely smooth.

\bl\label{l:kernel} If $\Delta x_F \in \rn$,  $\Delta x_F \perp c_F$
for all $F$, then $\{\Delta x_F\}\in {\rm Ker} D_{\{x^0_F\}} g$,
where $D_X g$ is the differential of $g$ at the point $X$. \el \bp
Since the kernel of the differential is a linear space, it suffices
to check this for the vectors $\{\Delta x_F\}$ in which only one
$\Delta x_{\widetilde{F}} \not=0$.
 Due to symmetry, we may assume that $c_{\widetilde{F}}=(\underbrace{1, \dots, 1}_{k}, \underbrace{0,\dots,0}_{n-k})$. The space orthogonal to $c_{\widetilde{F}}$ is then
 generated by the vectors $e_j$, $j>k$ and  $e_i-e_j$,  $1\le i <j \le k$.  Note now that the polytopes $Q^+$ and $Q^-$ built on the points $x_F^0$,  $F\not=\widetilde{F}$
 and $x_{\widetilde{F}}^0 \pm h e_j$, where $j>k$, are symmetric with respect to the symmetry $e_j \to -e_j$, so their volumes are the same.  On the other hand, the difference
 of their volumes in the first order is $2h D_{\{x^0_F\}}g (\{0, \dots, e_j,\dots, 0\})$, where $e_j$ stands in the position corresponding to $\widetilde{F} \in {\mathcal F}$.
 Thus,
$$
D_{\{x^0_F\}}g (\{0, \dots, e_j,\dots, 0\})=0.
$$
To prove the equality $D_{\{x^0_F\}}g (\{0, \dots, e_i-e_j,\dots, 0\})=0$, consider $Q'$ and $Q''$ built using the points $x_F=x_F^0$, $F \not=\widetilde{F}$ and  $x_{\widetilde{F}}=x_{\widetilde{F}}^0+h e_i$ or $x_{\widetilde{F}}=x_{\widetilde{F}}^0+h e_j$ respectively. They are also symmetric with respect to the symmetry $e_i \leftrightarrow e_j$ and the difference of their volumes in the first order equals
$hD_{\{x^0_F\}}g (\{0, \dots, e_i-e_j,\dots, 0\})$.

\ep

Below we shall also need the following elementary observation from
real analysis.
 \bl\label{l:real} Let $g(X)$ be a smooth function on
${\mathbb R}^N$, $X_0, X_1, X_2 \in {\mathbb R}^N$,  and
$$\|X_1-X_0\|, \|X_2-X_0\| \le \delta \to 0.$$
 Suppose that $X_1-X_2 \in
{\rm Ker} D_{X_0}g$. Then $|g(X_1)-g(X_2)| \le C \delta^2$.

\el
\bp Using the Taylor formula, we get
$$
g(X_j)=g(X_0)+\left(D_{X_0}g\right)(X_j-X_0)+O(\delta^2), \mbox{
where } j=1,2.
$$
Subtracting these two identities, we obtain
$$
g(X_1)-g(X_2)=\left(D_{X_0}g
\right)(X_1-X_2)+O(\delta^2)=O(\delta^2),
$$
because $\left(D_{X_0}g\right)(X_1-X_2)=0$. \ep

Let $P \subset \rn$ be a convex polytope.  For a face $F$ of $P$, we
 define its dual face  $F^*$ of $P^*$  by $F^*=\{y\in P^*: x\cdot y =1  \mbox{ for all    } x\in F \}$
 (see Chapter 3.4 in \cite{Gr}).

\bl\label{l:tochki}
 let $P$ be a convex polytope such that $0$ is in the interior of
 $P$. Let $P^*$ be its dual polytope. Chose some pair of dual faces
 $F$ and $F^*$ of $P$ and $P^*$ respectively and some points $x\in
 F$, $x^*\in F^*$ in the relative interiors of $F$ and $F^*$. Assume
 that $K$ is a convex body satisfying $(1-\delta) P \subset K
 \subset P$. Then there exists a pair of points $y\in \partial K$
 and $y^*\in \partial K^*$ such that $y \cdot y^*=1$ and $\|y-x\|,
 \|y^*-x^*\|\le C \delta$, where $C>0$ does not depend on $K$ or $\delta$,
 but may depend on $P, P^*, F, F^*, x$ and $x^*$.
 \el
\bp Since $x\cdot x^*=1>0$,  there exists a self-adjoint positive
definite linear operator $A$ such that $Ax=x^*$. This operator can
be chosen as follows: Let $L$ be a $2$-dimensional plane through the
origin containing both $x$ and $x^*$. $A$ will act identically on
$L^\perp$. To define its action on $L$, choose an orthogonal basis
$e_1, e_2$ in $L$ such that $e_1=x$ and put
$$
A\big|_{L}=\left(\begin{array}{cc}
     a & b \\
     b & a'
   \end{array}
   \right),
   $$
where $x^*=ae_1+be_2$ and $a'>0$ is chosen so large that $a a'>
b^2$.

 We will use below the following simple orthogonality
relations:
\begin{enumerate}
\item $x \perp l(F^*)$.
\item $x^* \perp l(F)$.
\item $l(F)\perp l(F^*)$.
\item $\left[ A^{-1} l(F^*) \right]^\perp ={\rm span}\left[x^*, A
l(F) \right]$ and  $\left[ A l(F) \right]^\perp ={\rm span}\left[x,
A^{-1} l(F^*) \right]$.
\item $(x^*)^\perp \cap {\rm span}(x, A^{-1} l(F^*))=A^{-1} l(F^*)$.
\end{enumerate}
$(1),(2)$ and $(3)$ follow directly from the definition of $F^*$
(see Chapter 3.4 in \cite{Gr}). Let us first prove (4). Since $l(F)
\perp l(F^*)$ and $A$ is self-adjoint, we also have $A l(F) \perp
A^{-1}l(F^*)$. Also, since $x \perp l(F^*)$, we have $x^*=Ax \perp
A^{-1}l(F^*)$. Thus ${\rm span}(x^*, A l(F))\subset \left[ A^{-1}
l(F^*)\right]^\perp$. On the other hand, $x\not \in l(F)$, so $x^*=
Ax \not\in A l(F)$ and
$$
{\rm dim} \left({\rm span}(x^*, Al(F))\right)=1+ {\rm dim} F=n -
{\rm dim} F^*=n-{\rm dim} A^{-1} l(F^*),
$$
so $A^{-1}l(F^*)^{\perp}$ can not be wider than ${\rm span}(x^*,
Al(F))$. Similarly,
$$\left[ A l(F) \right]^\perp ={\rm
span}\left[x, A^{-1} l(F^*) \right].$$

To prove (5), we first note that $A^{-1} l(F^*)\perp x^*$ (see (4)).
Since $x^* \cdot x=1 \not =0$, $(x^*)^\perp \cap {\rm span}(x,
A^{-1} l(F^*))$ is a subspace of codimension $1$ in ${\rm span}(x,
A^{-1} l(F^*))$, so it cannot be wider than $A^{-1}l(F^*)$.

Let $\widetilde{K}=K\cap {\rm span}(x, A^{-1} l(F^*))$ and let $y
\in \widetilde{K}$ maximize $y \cdot x^*$. Then $y \in
\partial K$ and a tangent plane to $K$ at $y$ contains an affine
plane parallel to $(x^*)^\perp \cap {\rm span}(x, A^{-1} l(F^*))=
A^{-1} l(F^*)$. Therefore, there exists $y^* \in \partial K^* \cap
\left[A^{-1} l(F^*) \right]^\perp=\partial K^*\cap{\rm span} (x^*, A
l(F))$ such that $y\cdot y^*=1$.

Now let $y=\alpha x +h$ and $y^*=\alpha^*x^*+h^*$,  where  $h\in
A^{-1}l(F^*)$ and $h^*\in A l(F)$. Note that $y \cdot x^* = \alpha$,
so by maximality of $y$,
$$
\alpha=(y,x^*) > (0, x^*)=0.
$$
Also $y\cdot y^*=\alpha \alpha^*=1$, so $\alpha^* >0$. Let  $\rho >
0$ be such that $B(x,\rho)\cap \af (F) \subset F$  and  $B(x^*,
\rho)\cap \af (F^*) \subset F^*$ where $B(z,t)$ is the Euclidean
ball of radius $t$ centered at $z$. Since $y\in
\partial K$ and
$$
K^* \supset P^* \supset F^* \ni x^* + \frac{\rho A h}{\|Ah\|},
$$
we  have
$$
1 \ge y \cdot \left(x^* + \frac{\rho A h}{\|Ah\|}\right) = \alpha +
\frac{\rho A h \cdot h}{\|Ah\|} \ge \alpha + \rho' \|h\|, \mbox{
where  } \rho'=\frac{\rho}{\|A\|\|A^{-1}\|}.
$$
Since $y^* \in \partial K^*$ and
$$
K\supset (1-\delta) P \supset (1-\delta) F \ni (1-\delta)\left[x +
\frac{\rho A^{-1} h^*}{\|A^{-1}h^*\|}\right],
$$
we have
$$\left( (1-\delta)\left[x + \frac{\rho A^{-1}
h^*}{\|A^{-1}h^*\|}\right]\right) \cdot y^* \le 1
$$
 and
$$
(1-\delta)^{-1} \ge \left[x + \frac{\rho A^{-1}
h^*}{\|A^{-1}h^*\|}\right] \cdot y^* = \alpha^* +\frac{\rho A^{-1}
h^* \cdot h^*}{\|A^{-1}h^*\|} \ge \alpha^*+\rho' \|h^*\|.
$$
Thus $\alpha  \le 1$ and $\alpha^*\le 1/(1-\delta)$, which, together
with $\alpha \alpha^*=1$, gives $\alpha \ge 1-\delta$ and
$\alpha^*\ge 1$. Hence $\rho'\|h\|\le \delta$, $\rho'\|h^*\| \le
\frac{1}{1-\delta}-1$ and, thereby, $\|y-x\|, \|y^*-x^*\|\le C
\delta$. \ep

 Now define $c_F^*=\frac{1}{n-{\rm dim} F}c_F$. Choose positive
 numbers $\alpha_F$ and $\alpha^*_F$ satisfying
 $\alpha_F\alpha^*_F=1$ and put $y_F=\alpha_F c_F$, $y_F^*=\alpha_F^*
 c_F^*$.

 Let $Q=\cup_{\mathbb F} S_{\mathbb F}(Q)$ and $Q'=\cup_{\mathbb F} S_{\mathbb
 F}(Q')$, where
$$
S_{\mathbb F}(Q)={\rm conv}(0, y_{F_0}, y_{F_1}, \dots, y_{F_{n-1}})
\mbox{ and } S_{\mathbb F}(Q')={\rm conv}(0, y^*_{F_0}, y^*_{F_1},
\dots, y^*_{F_{n-1}})
$$
and ${\mathbb F}$ runs over all flags ${\mathbb  F}=\{F_0, \dots,
F_{n-1}\}$ of faces of $B_\infty^n$.

\bl\label{l:c}
$$
\vol_n(Q)\vol_n(Q') \ge \Prod(B_\infty^n).
$$
\el

\bp For every flag ${\mathbb  F}=\{F_0, \dots, F_{n-1}\}$,
$$
\vol_n(S_{\mathbb  F}(Q))=\vol_n (S_{\mathbb
F}(B_\infty^n))\prod\limits_{j=0}^{n-1} \alpha_{F_j}, \mbox{ where }
S_{\mathbb F}(B_\infty^n)={\rm conv}(0, c_{F_0}, c_{F_1},\dots,
c_{F_{n-1}}),
$$
and
$$
\vol_n(S_{\mathbb  F}(Q'))=\vol_n (S_{\mathbb
F}(B_1^n))\prod\limits_{j=0}^{n-1} \alpha_{F_j}^*, \mbox{ where }
S_{\mathbb F}(B_1^n)={\rm conv}(0, c^*_{F_0}, c^*_{F_1},\dots,
c^*_{F_{n-1}}).
$$
Hence
$$
\vol_n(S_{\mathbb  F}(Q)) \vol_n(S_{\mathbb  F}(Q')) =  \vol_n
(S_{\mathbb  F}(B_\infty^n))\vol_n (S_{\mathbb  F}(B_1^n)).
$$
The factors on the right hand side do not depend on the flag
${\mathbb F}$. Thus,
$$
\vol_n(Q)\vol_n(Q')=\sum\limits_{\mathbb  F}\vol_n(S_{\mathbb
F}(Q))\sum\limits_{\mathbb  F}\vol_n(S_{\mathbb F}(Q'))
$$
$$
\ge \left(\sum\limits_{\mathbb  F}\sqrt{\vol_n(S_{\mathbb
F}(Q))\vol_n(S_{\mathbb F}(Q'))} \right)^2
=\left(\sum\limits_{\mathbb F}\sqrt{\vol_n(S_{\mathbb
F}(B_\infty^n))\vol_n(S_{\mathbb F}(B_1^n))} \right)^2
$$
$$=\sum\limits_{\mathbb  F}\vol_n(S_{\mathbb
F}(B_\infty^n))\sum\limits_{\mathbb  F}\vol_n(S_{\mathbb F}(B_1^n))=
\vol_n(B_\infty^n)\vol_n(B_1^n)=\Prod(B_\infty^n).
$$
 \ep

\section{Lower stationarity of $B_\infty^n$}

Now apply Lemma \ref{l:tochki} to $B_\infty^n$ and $B_1^n$ and the
points $c_F \in F$ and $c^*_F=\frac{1}{n-{\rm dim} F} c_F \in F^*$,
where $F^*$ is the face of $B_1^n$ dual to $F$. Since in this case
we can choose $A$ (in the proof of Lemma \ref{l:tochki}) to be a
pure homothety with coefficient $\frac{1}{n-{\rm dim} F}$, we get
points
$$
x_F=\alpha_Fc_F+h_F \mbox{  and  } x_F^*=\alpha_F^*c_F^*+h_F^*
\mbox{  satisfying } x_F\in \partial K, x_F^* \in \partial K^*,
$$
where $\alpha_F \alpha^*_F=1$, $h_F \in l(F^*)$, $h^*_F\in l(F)$ and
$|\alpha_F-1|, |\alpha^*_F-1|, \|h_F\|, \|h^*_F\| \le C\delta$.

Since $\pm e_j \in \partial K$ and $\pm e_j \in \partial K^*$, we
can choose $x_F=y_F=x^*_F=y_F^*=c_F=c^*_F$ when ${\rm dim} F=n-1$.

Put $y_F=\alpha_F c_F$ and $y^*_F=\alpha_F^* c_F^*$, and consider
the polytopes
$$
P=\cup_{\mathbb F} {\rm conv}(0, x_{F_0},\dots,  x_{F_{n-1}}) \mbox{
and } P'=\cup_{\mathbb F} {\rm conv}(0, x_{F_0}^*,\dots,
x^*_{F_{n-1}}),
$$
$$
Q=\cup_{\mathbb F} {\rm conv}(0, y_{F_0},\dots,  y_{F_{n-1}}) \mbox{
and } Q'=\cup_{\mathbb F} {\rm conv}(0, y^*_{F_0},\dots,
y^*_{F_{n-1}}).
$$
Note that $x_F-y_F=h_F$, $x^*_F-y_F^*=h_F^*$ and $h_F, h^*_F \perp
c_F$.

Thus by Lemmata \ref{l:kernel}, \ref{l:real}.
$$
|\vol_n(P)-\vol_n(Q)| \le C\delta^2 \mbox{ and }
|\vol_n(P')-\vol_n(Q')| \le C\delta^2,
$$
whence
$$
\vol_n(P)\vol_n(P')\ge \vol_n(Q)\vol_n(Q')-C\delta^2 \ge
\Prod(B_\infty^n)-C\delta^2,
$$
where the last inequality follows from Lemma \ref{l:c}.

Since $K\supset P$ and $K^* \supset P'$, it remains to show that for
some $c>0$, either $K\not\subset (1+c\delta)P$, or $K^*\not\subset
(1+c\delta)P'$. Then, by Lemma \ref{l:small}, ether $\vol_n(K)\ge
\vol_n(P)+c'\delta$, or $\vol_n(K^*) \ge \vol_n(P')+c'\delta$. This
yields
$$
\Prod(K) \ge \Prod(B_\infty^n)+c''\delta-C\delta^2>
\Prod(B_\infty^n),
$$
provided that $\delta>0$ is small enough.

\section{The conclusion of the proof}

Note that at least one of the coordinates of one of the
$x_{\widetilde{F}}$ with ${\rm dim} \widetilde{F} =0$ is at most
$1-\delta$. Indeed, assume that all coordinates are greater then
$(1-\delta')$ in absolute value with some $\delta'<\delta$. Define
$D=\mbox{conv}\{x_{F}: F\in {\mathcal F}, {\rm dim}F =0\} \subset
K$. Let $z \in D^*$. Choose $F$ so that $(x_F)_j z_j \ge 0$ for
all $j=1, \dots, n$. Then
$$
1 \ge x_F \cdot z\ge (1-\delta')\sum\limits_j |z_j|.
$$
Thus $D^* \subset (1-\delta')^{-1} B_1^n$ and $D\supset
(1-\delta')B_\infty^n$, contradicting the minimality of $\delta$.

Due to symmetry, we may assume without loss of generality that
$\widetilde{F}=\{(1,\dots,1)\}$ and that $(x_{\widetilde{F}})_1 \le
1-\delta$. Assume that $K\subset (1+c \delta)P$. Consider the point
$\tilde{x}=(1-\delta, c'\delta, \dots, c' \delta)$ where
$c'=1/(n-\frac{5}{4})$. Then $\tilde{x} \in (1-c''\delta)P^*$, where
$c''=1/(4n-5)$. Indeed, it is enough to check that $\tilde{x} \cdot
x_F \le 1-c''\delta$ for all vertices $x_F$ of $P$. If
$F\not=\{(1,\dots, 1)\}$, then all coordinates of $x_F$ do not
exceed $1$ and at least one does not exceed $1/2$. Thus, if $\delta$
is small enough, we get
$$
\tilde{x}\cdot x_F \le (1-\delta)+(n-2)c'\delta+\frac{c'\delta}{2} =
1-\delta+(n-\tfrac{3}{2})c'\delta =1-c''\delta.
$$
If $F=\{(1,\dots,1)\}$, then
$$
\tilde{x}\cdot x_F \le
(1-\delta)^2+(n-1)c'\delta=1-2\delta+\frac{n-1}{n-\tfrac{5}{4}}
\delta+\delta^2 \le 1-2\delta+\frac{4}{3}\delta+\delta^2\le
1-c''\delta,
$$
provided that $\delta>0$ is small enough.

Therefore if $c < c''$, we get $\tilde{x}\in
\frac{1}{1+c\delta}P^*\subset K^*$.

Now note that for every $x\in P'$, we have
$$
|x_1|+(1-C'\delta)\sum\limits_{j \ge 2}|x_j| \le 1,
$$
provided $C'$ is chosen large enough. Indeed, again it is enough to
check this for the vertices $x^*_F$ of $P'$. If $c_F \not=(\pm
1,0,\dots,0)$ we have $\sum\limits_{j \ge 2}|(x^*_F)_j|\ge 1/3$, so
$$
|(x^*_F)_1|+(1-C'\delta)\sum\limits_{j \ge 2}|(x^*_F)_j| \le
\sum\limits_{j\ge 1}|(x^*_F)_j|-C'\delta\sum\limits_{j \ge 2}|(x^*_F)_j| \le 1+n
C\delta-\frac{C'\delta}{3} \le 1,
$$
provided that $C' \ge 3n C$, where $C$ is the constant such that
$\|x_F^*-c_F^*\|\le C\delta$. If $c_F=(\pm 1, 0\dots,0)$, then
$x_F=\pm e_1$ and the inequality is trivial.

Now it remains to note that
$$
|\tilde{x}_1|+(1-C'\delta)\sum\limits_{j \ge 2}|\tilde{x}_j|
=1-\delta+(1-C'\delta)(n-1)c'\delta=1+c''\delta-C'(n-1)c'\delta^2 >
1+c\delta,
$$
provided that $c<c''/2$ and $\delta$ is small enough, whence
$\tilde{x} \not\in (1+c\delta)P'$.

\end{document}